\documentclass{amsart}

\usepackage{verbatim}
\usepackage{enumitem}
\usepackage{xcolor}
\usepackage{tikz}
\usepackage{tikz-cd}
\usepackage{arydshln}
\usepackage{caption}
\usepackage{subcaption}
\usepackage{graphicx}
\usepackage{ulem}

\tikzset{
	symbol/.style={
		draw=none,
		every to/.append style={
			edge node={node [sloped, allow upside down, auto=false]{$#1$}}}
	}
}

\usepackage{amsmath}
\usepackage{bbm}
\usepackage{hyperref}
\usepackage{times,amsfonts,amsmath,amstext,amsbsy,amssymb,
amsopn,amsthm,upref,eucal}
\usepackage[T1]{fontenc}

\newtheorem{theorem}{Theorem}[section]

\newtheorem{corollary}[theorem]{Corollary}

\theoremstyle{theorem}
\newtheorem{definition}[theorem]{Definition}

\theoremstyle{remark}
\newtheorem{remark}[theorem]{Remark}

\numberwithin{equation}{section}

\begin{document}

\title[Chaos in large genus surfaces]{Chaos in large genus surfaces}



\author{Francisco Arana--Herrera}

\email{fa50@rice.edu}

\address{Department of Mathematics, Rice University,
Herman Brown Hall for Mathematical Sciences
6100 Main Street
Houston, TX 77005.}



\date{}

\begin{abstract}
    Geodesic flows on closed hyperbolic surfaces are a quintessential example of chaotic dynamics, i.e., systems whose long term behavior is very sensitive to initial conditions. The speed of such chaos is controlled by the spectral gap of the Laplace--Beltrami operator of the underlying hyperbolic surface. In this paper we give an overview of recent breakthroughs of Anantharaman and Monk showing that large genus closed hyperbolic surfaces have optimal spectral gap in a probabilistic sense. On the way we introduce and discuss the foundational works of many authors, from Selberg to Mirzakhani, that play a crucial role in the \textit{tour de force} proof of Anantharaman and Monk.
\end{abstract}

\maketitle


\thispagestyle{empty}

\tableofcontents

\section{Introduction}\label{sec:intro}

\subsection{Chaos.}\label{subsec:chaos} 
The main purpose of classical mechanics is to study physical systems that evolve under deterministic laws of motion: knowing the current state of a system completely determines its future. Although in theory this appears to be a rather reasonable goal, under further consideration one runs into a practical issue of the utmost importance: If real life measurements are only an approximation of the underlying physical quantities, how can one even claim to know the current state of a system?

At first one might think that this is not really an issue: Good enough approximations of the initial conditions of a system should lead to reasonable approximations of the outcome. And although this is certainly true if the equations of motion are what one calls \textit{well-posed}, the dependence of the outcomes on the initial conditions of the system can be quite dramatic. This type of behavior is usually referred to as \textit{chaos}.

\begin{figure}[ht]
\centering
\includegraphics[scale=.18]{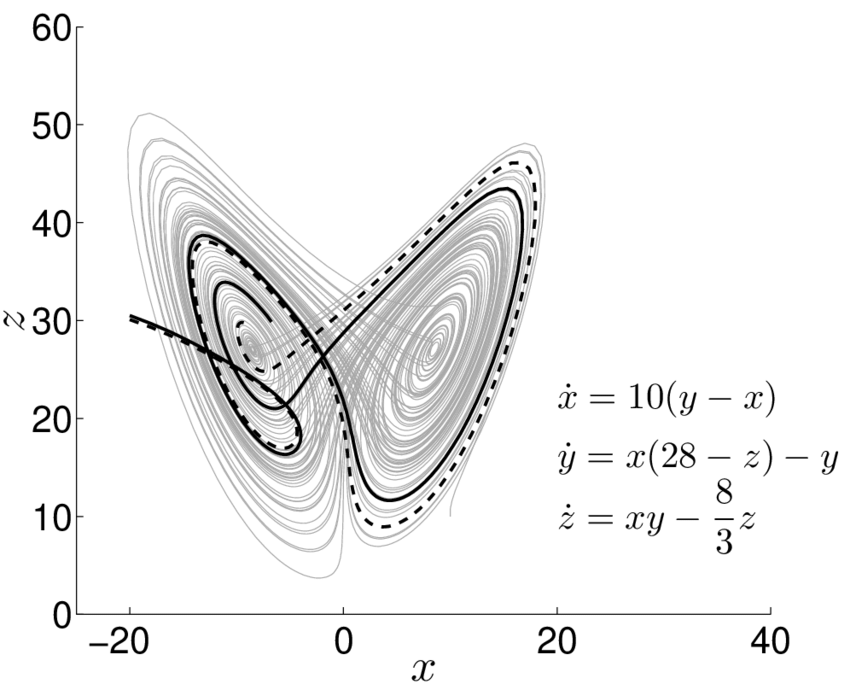}
\caption{An instance of the Lorenz attractor \cite{lorenz}.}
\label{fig:lorenz}
\end{figure}

And to be emphatic, this is certainly an issue we face in our everyday lives: How many times have you trusted the weather forecast and left home without an umbrella only to find yourself caught unprepared in the middle of a thunderstorm only a few hours later? The earth's atmosphere is a \textit{chaotic system}, but a rather complicated one. Although it might seem counterintuitive at first, it turns out that systems with few degrees of freedom can also be chaotic. Perhaps the most famous example of such a system is given by the \textit{Lorenz attractor}; see Figure \ref{fig:lorenz} for a particular instance of this example.

In this setting, let us formulate a few questions that will guide our discussion:
\begin{itemize}
    \item \textbf{What mechanisms give rise to chaos?}
    \item \textbf{Can we quantify the chaoticity of a system?}
    \item \textbf{What controls the speed of chaos?}
    \item \textbf{Do random systems have a tendency to be more chaotic?}
\end{itemize}

\subsection{Hyperbolic surfaces.} Certainly, studying these questions in full generality would be an extremely ambitious goal. So, as we do many times in mathematics, we will aim at understanding these questions by restricting to a class of concrete toy models: \textit{geodesic flows on closed hyperbolic surfaces}. The hope here is that, once we understand these toy models well enough, it should be possible to extrapolate general principles to more complicated systems. But let us not jump too far ahead of ourselves. In this survey we will explore the following answers to the questions above for closed hyperbolic surfaces:\medskip

\begin{itemize}
    \item \textbf{What mechanisms give rise to chaos?} \newline
    The \textit{negative curvature} of a closed hyperbolic surface forces the orbits of its geodesic flow to drift apart quickly, giving rise to chaotic behavior. \medskip
    \item \textbf{Can we quantify the chaoticity of a system?} \newline
    The geodesic flow on a closed hyperbolic surface is \textit{exponentially mixing}, i.e., the future and past become independent at exponential speed.
    \medskip 
    \item \textbf{What controls the speed of chaos?} \newline
    The speed of mixing of the geodesic flow on a closed hyperbolic surface is controlled by the smallest non-zero eigenvalue of its \textit{Laplace--Beltrami operator}. \medskip
    \item \textbf{Do random systems have a tendency to be more chaotic?} \newline
    Random large genus closed hyperbolic surfaces have \textit{optimal spectral gap}. \medskip
\end{itemize}

\subsection{Organization.} Our main aim in this paper is to discuss the recent breakthroughs of Anantharaman and Monk \cite{gap1,monksur,gap2} answering the last question above. Although this work is quite technical, worry not, as we will also spend  time answering the first three questions to set up notation and introduce the subject to the unfamiliar reader. 

In \S2 we discuss the basics of hyperbolic surfaces. In \S3 we introduce the notion of chaos we will be interested in: exponential mixing. In \S4, following Ratner, we discuss how the spectrum of the Laplace--Beltrami operator controls the mixing speed of the geodesic flow of a closed hyperbolic surface. In \S5 we introduce the Selberg trace formula, one of the key ingredients in the work of Anantharaman and Monk. In $\S6$ we discuss the notion of Weil--Petersson random hyperbolic surface and Mirzakhani's foundational work on the subject; this is another key ingredient in the work of Anantharaman and Monk. Finally, in \S7 we give a complete statement of the main result and describe the new key techniques that make the proof of Anantharaman and Monk so revolutionary.

\subsection{Acknowledgments.} This expository article was prepared to accompany a talk given by the author at the Joint Mathematics Meetings in 2026. The author would like to thank Bianca Viray and Daniel Erman for the invitation to give this talk. The author would also like to thank Giovanni Forni for many enlightening conversations.

\section{Hyperbolic surfaces}

\subsection{Hyperbolic surfaces.} 
\begin{figure}[ht]
\centering
\includegraphics[scale=.35]{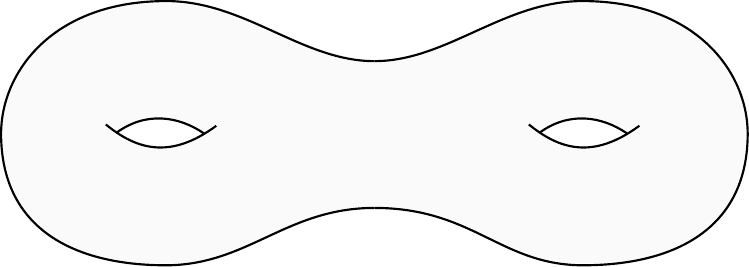}
\caption{An orientable surface of genus $2$.}
\label{fig:g2}
\end{figure}

Hyperbolic surfaces are what we call a \textit{homogeneous space}, i.e., a space whose local geometry is the same everywhere. In the case of hyperbolic surfaces, the local geometry is modeled on the hyperbolic plane. Roughly speaking, near every point the surface looks like a saddle and the bending of such saddles is the same across the whole surface. Closed orientable surfaces are completely classified by their \textit{genus}, i.e., the number of holes; see Figure \ref{fig:g2} for an example. In the case of closed hyperbolic surfaces, the genus $g$ is related to the area $A$ by the Gauss--Bonnet theorem:
\[
A = 4\pi(g-1).
\]

\begin{remark}
    For the rest of this paper we exclusively consider orientable surfaces and just refer to them as \textit{surfaces}, without highlighting their orientability.
\end{remark}

But we are jumping too far ahead of ourselves! Let us spend some time making some of the basic notions in hyperbolic geometry more precise.

\subsection{The hyperbolic plane.} 
\begin{figure}[ht]
\centering
\includegraphics[scale=.75]{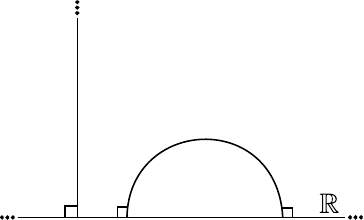}
\caption{Examples of geodesics on the hyperbolic plane.}
\label{fig:h2}
\end{figure}

The \textit{hyperbolic plane} can be defined as the upper half space
\[
\mathbb{H} := \{z \in \mathbb{C} \ | \ \Im(z) > 0 \}
\]
endowed with the Riemannian metric given by
\[
g := \frac{dx^2 + dy^2}{y^2};
\]
this is nothing more to say than the norm of a tangent vector $v := (v_1,v_2) \in \mathbb{R}^2$ at a point $z = x+iy \in \mathbb{H}$ is given by
\[
\|v\|_z := \frac{\sqrt{v_1^2 + v_2^2}}{y}.
\]
Endowed with this geometry, the \textit{geodesics}, i.e., the length minimizing curves, of $\mathbb{H}^2$ are the vertical lines and the half circles of the upper half space with center along the real axis; see Figure \ref{fig:h2} for some examples. As one can measure lengths on $\mathbb{H}$, one can also measure areas. More explicitly, the corresponding infinitesimal area element is given by
\[
dA =: \frac{dx dy}{y^2}.
\]

For any of this to make sense, one needs to make sure that the geometry near every point of $\mathbb{H}$ is the same. Said another way, one needs to guarantee that the group of isometries of $\mathbb{H}$ acts \textit{transitively} on its unit tangent bundle, i.e., has a single orbit. And indeed that is the case. Furthermore, one can show that, up to finite index, the group of orientation preserving isometries of $\mathbb{H}$ is given by
\[
\mathrm{SL}(2,\mathbb{R}) := \left\lbrace \left( \begin{array}{c c}
a & b \\ c & d
\end{array}\right) \in \mathrm{Mat}_{2\times2}(\mathbb{R})\ \bigg\vert \ a,b,c,d \in \mathbb{R}, \thinspace ab-cd = 1 \right\rbrace
\]
acting by \textit{Möbius transformations} on $\mathbb{H}$ in the following way:
\[
\left(\begin{array}{c c}
a & b \\
c & d
\end{array}\right) \cdot z := \frac{az+b}{cz+d}, \quad z \in \mathbb{H}^2.
\]

\subsection{Negative curvature.} 
\begin{figure}[ht]
		\centering
		\begin{subfigure}{.45\textwidth}
			\centering
			\includegraphics[scale=.35]{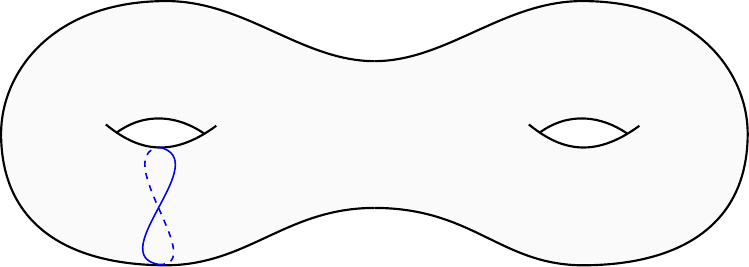}
		\end{subfigure}
		\begin{subfigure}{.45\textwidth}
			\centering
			\includegraphics[scale=.35]{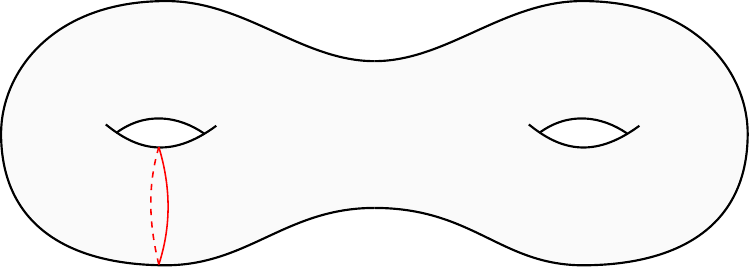}
		\end{subfigure}
        \caption{Tightening a loop (left) to its geodesic representative (right).}
        \label{fig:tight}
	\end{figure}

Now notice that the geodesics represented in Figure \ref{fig:h2} seem to drift apart rather quickly; compare this to the more familiar Euclidean setting where straight lines that do not intersect always stay at the same distance.
Indeed, one can show that geodesics in $\mathbb{H}$ drift away at exponential speed. This is one of the fundamental manifestations of \textit{negative curvature}, i.e., the saddle-like behavior we described at the beginning of this section, and will play a fundamental role in the rest of our discussion. An important consequence of negative curvature is the fact that every loop on a closed hyperbolic surface can be freely homotoped into a unique geodesic loop. Such a geodesic representative can be obtained by tightening the original loop as to minimize its length; see Figure \ref{fig:tight}.

\begin{remark}
    To the concerned reader we point out that, although the closed hyperbolic surfaces we draw, e.g., the ones in Figure \ref{fig:tight}, are certainly not geometrically accurate, in the sense that they do not look negatively curved everywhere, they will help us explain some of the concepts that would otherwise be harder to present.
\end{remark}

\subsection{The geodesic flow.} Given a closed hyperbolic surface $X$ we denote by $T^1X$ its corresponding unit tangent bundle. The \textit{geodesic flow} on this bundle,
\[
\Phi := \{\phi_t \colon T^1 X \to T^1 X\}_{t \in \mathbb{R}},
\]
is the flow defined by following the geodesic spanned by a given unit tangent vector for $t$ units of length. The orbits of this flow, or rather their projections to $X$, are precisely the geodesics of the surface. These flows constitute the main dynamical systems whose chaotic behavior we will be interested in understanding.

\section{Exponential mixing}

\subsection{Random processes.} To introduce the notion of chaos we will be interested in, let us first explore an analogous notion in the context of probability theory. Recall that, given a probability space $\Omega$ with a probability measure $\mathbb{P}$, two measurable events $A,B \subseteq \Omega$ are said to be \textit{independent} if they satisfy
\[
\mathbb{P}(A \cap B) = \mathbb{P}(A) \cdot \mathbb{P}(B).
\]

Now, in the same context, we could have a random process, i.e., a sequence of random variables $(Y_t)_{t \in \mathbb{R}}$, say, for simplicity, all distributed like $Y$, and ask whether the present and future become independent as we let the time grow larger and larger. More formally this condition can be stated as follows: for every pair of measurable sets $A,B \subseteq \mathbb{R}$,
\[
\lim_{t \to \infty} \mathbb{P}(Y_0 \in A, \ Y_t \in B) = \mathbb{P}(Y\in A) \cdot \mathbb{P}(Y \in B).
\]

\subsection{Mixing.} The condition above is the one we want to capture in the setting of geodesic flows on closed hyperbolic surfaces. Let $X$ be a closed hyperbolic surface and $\Phi$ be the corresponding geodesic flow on its unit tangent bundle $T^1X$. To proceed in analogy with the discussion above, we need a probability measure on our phase space, i.e., on $T^1X$. Although many choices are possible, we will stick with the most natural one, the \textit{Liouville measure}, that is, the probability measure $\mu$ corresponding to the area induced by the hyperbolic metric on $X$ and which is uniform with respect to angle on fibers of $T^1X$.

In this setting we say that the geodesic flow $\Phi$ on $T^1X$ is \textit{mixing} if for every pair of measurable sets $A,B \subseteq T^1 X$,
\[
\lim_{t \to \infty} \mu(\{v \in T^1X \ | \ v \in A, \ \phi_t(v) \in B)\}) = \mu(A) \cdot \mu(B).
\]
For later reference, let us note that this condition can also be written as
\[
\lim_{t \to \infty} \int_{T^1X} \mathbbm{1}_A(v) \thinspace \mathbbm{1}_B(\phi_t(v)) \thinspace d\mu(v) = \int_{T^1X} \mathbbm{1}_A(v) \thinspace d\mu(v) \cdot \int_{T^1X} \mathbbm{1}_B(v) \thinspace d\mu(v).
\]
We interpret this condition as chaos because
it guarantees that small changes in the initial condition, i.e., even if we allow it to vary in a tiny positive measure set, produce wildly different outcomes in the long term; indeed, every possible outcome with the correct proportion given by the underlying Liouville measure.

\subsection{Exponential mixing.} Now, the asymptotic independence captured by the notion of mixing could be attained at different speeds. This question is a quite delicate one, and, indeed, it is hard to make sense of it unless one restricts to sufficiently regular observables. To this end, we will consider continuously differentiable functions $f,g \in C^1(T^1X)$ and denote their respective $C^1$ norms by $\|f\|_{C^1}$ and $\|g\|_{C^1}$. In this setting we say that the geodesic flow $\Phi$ on $T^1X$ is \textit{exponentially mixing}, or chaotic with exponential speed, if there exists constants $C > 0$ and $\kappa > 0$ such that for such every pair $f,g \in C^1(T^1X)$,
\[
\bigg\vert  \int_{T^1X} f(v) \thinspace g(\phi_t(v)) \thinspace d\mu(v) - \int_{T^1X} f(v) \thinspace d\mu(v) \cdot \int_{T^1X} g(v) \thinspace d\mu(v) \bigg\vert \leq C \|f\|_{C^1}\|g\|_{C^1} e^{-\kappa t}.
\]

In this setting we have the following result originally due to Moore:

\begin{theorem}[\cite{Moore}]
    \label{theo:moore}
    The geodesic flow $\Phi$ on the unit tangent bundle $T^1X$ of any closed hyperbolic surface $X$ is exponentially mixing.
\end{theorem}

Moore's original proof uses the fact that the geodesic flow $\Phi$ is the diagonal part of a larger $\mathrm{SL}(2,\mathbb{R})$ action, which provides a rich algebraic structure to study its properties. We highlight that, although this perspective will be very important for the discussion in the next section, exponential mixing is a much more prevalent phenomena, even beyond the homogeneous setting considered here. For instance, Dolgopyat \cite{dol} proved exponential mixing for geodesic flows on closed surfaces of variable negative curvature. 

\begin{remark}
Although we do not elaborate more on it here, we note that exponential mixing in general is an extremely powerful tool in applications. For instance, it was used by Kontorovich and Oh \cite{apol} to solve important sphere packing counting problems.
\end{remark}

\section{The Laplace--Beltrami operator}

\subsection{The hyperbolic plane.} For reasons that will become clear later, we want to define an analogue of the Euclidean Laplacian on $\mathbb{H}$. Recall that the gradient of a real valued function specifies the direction and rate of maximal increase. Recall also that the divergence of a vector field represents the infinitesimal change in volume induced by the corresponding flow. The Laplacian of a function is nothing more than the the divergence of its gradient vector field; it is also known as the \textit{Laplace--Beltrami operator}.

Just as this makes sense in Euclidean space, it also makes sense in the hyperbolic plane; one just needs to be careful to use the corresponding metric and induced volume form when computing the gradient and divergence. More explicitely, in the model considered above, the Laplace--Beltrami operator on $\mathbb{H}$ is the differential operator given by
\[
\Delta := y^2 \left(\frac{\partial^2}{\partial x^2} + \frac{\partial^2}{\partial y^2} \right).
\]

Denote by $L^2(\mathbb{H})$ the space of square-integrable functions on the hyperbolic plane with respect to the corresponding volume form. Then, $\Delta$ is an unbounded self-adjoint operator on $L^2(\mathbb{H})$ whose \textit{spectrum} is the interval $(-\infty,-1/4]$; the reader unfamiliar with the notion of spectrum can think of it as a generalization of eigenvalues to infinite dimensional vector spaces. The number $1/4$ will play a crucial role in the rest of this paper. 

\begin{remark}
    We refer the reader to \cite{plane_gap} for a fantastic exposition of these and more general results on the spectral theory of Laplace--Beltrami operators.
\end{remark}

\subsection{Closed hyperbolic surfaces.} Because the Laplace--Beltrami operator of $\mathbb{H}$ is defined only in terms of its local geometry, it is invariant under its group of isometries. In particular, one can define the Laplace--Beltrami operator on the space $L^2(X)$ of square-integrable functions of an arbitrary closed hyperbolic surface $X$ via its local identifications with $\mathbb{H}$. We will also denote this operator by $\Delta$, making it clear from the context when it corresponds to the operator on $L^2(X)$ rather than the one on $L^2(\mathbb{H})$.

The Laplace--Beltrami operator $\Delta$ on a closed hyperbolic surface $X$ is thus an unbounded self-adjoint operator on $L^2(X)$. Furthermore, $\Delta$ has compact resolvent and thus can be diagonalized in the sense that $L^2(X)$ admits a complete orthonormal system of eigenfunctions. More precisely, there exist real numbers and smooth functions
\begin{gather*}
    0 = \lambda_0 > \lambda_1 \geq \lambda_2 \geq \cdots \text{ with } \lim_{n \to \infty} \lambda_n = -\infty, \ \quad f_0,f_1,f_2,\dots \in C^\infty(X),
\end{gather*}
such that the following conditions hold:
\begin{enumerate}
    \item $L^2(X) = \overline{\mathrm{span}(\{f_n\}_{n \in \mathbb{N}})}$,
    \item $\Delta f_n = \lambda_n f_n \ \text{for all } n \in \mathbb{N}$.
\end{enumerate}

The first non-zero eigenvalue $\lambda_1 < 0$ of the Laplace--Beltrami operator $\Delta$ on $X$ will be of particular importance to us. We will denote this eigenvalue by $\lambda_1(X)$, making the dependence on $X$ explicit. We will also refer to $-\lambda_1(X) > 0$ as the \textit{spectral gap} of $X$.

At this stage one might ask if there is any relation between the spectrum of the Laplacian on a closed hyperbolic surface $X$ and the spectrum of the Laplacian on the hyperbolic plane $\mathbb{H}$. Such connections are not automatic because the square-integrability condition on the functions considered does not allow one to lift them from $X$ to $\mathbb{H}$. Still, it is possible to obtain the following topological bound on the spectral gap of any hyperbolic surface:

\begin{theorem} [\cite{comp2,comp1}]
    \label{theo:opt}
    For every $g \geq 2$ there exists $\epsilon(g) < 0$ with $\epsilon(g) \nearrow 0$ as $g \to \infty$ such that if $X$ is a closed hyperbolic surface of genus $g$ and $\lambda_1(X) < 0$ is the first non-zero eigenvalue of its Laplace--Beltrami operator $\Delta$, then
    \[
    \lambda_1(X) \geq -1/4 + \epsilon(g).
    \]
\end{theorem}

\begin{remark}
The interested reader can consult \cite{buser} for many other foundational results on the spectral theory of the Laplace--Beltrami operator on closed hyperbolic surfaces.
\end{remark}

\subsection{Speed of mixing.}
Although Theorem \ref{theo:moore} guarantees the geodesic flow on the unit tangent bundle of any closed hyperbolic surface is exponentially mixing, it does not provide explicit control on the speed of mixing. The following result of Ratner solves this issue by relating such speed to the spectral gap of the corresponding Laplace--Beltrami operator:

\begin{theorem}[\cite{ratner}]
    \label{theo:mix_rate}
    There exists a universal constant $C > 0$ with the following property. Let $X$ be a closed hyperbolic surface, $\lambda_1(X) < 0$ be the first non-zero eigenvalue of its Laplace Beltrami operator $\Delta$, and
    \[
    \sigma :=  1 - \sqrt{\max\{1 + 4 \lambda_1(X),0\}}.
    \]
    Then, for every pair of regular observables $f,g \in C^1(T^1X)$,
    \[
    \bigg\vert  \int_{T^1X} f(v)  g(\phi_t(v)) d\mu(v) - \int_{T^1X} f(v) d\mu(v) \cdot \int_{T^1X} g(v)  d\mu(v) \bigg\vert \leq C \|f\|_{C^1}\|g\|_{C^1} te^{-\sigma t}.
    \]
\end{theorem}

\begin{remark}
    It is also possible to control the spectral gap of a closed hyperbolic surface in terms of the exponential mixing rate of its geodesic flow; see for instance \cite{revrat}.
\end{remark}

Ratner actually proves a more general result about the decay of matrix coefficients of irreducible $\mathrm{SL}(2,\mathbb{R})$ representations in terms of the spectrum of the corresponding Casimir operator. Recall that, as highlighted by Moore's work, the geodesic flow $\Phi$ of a closed hyperbolic surface is the diagonal part of a larger $\mathrm{SL}(2,\mathbb{R})$ action. Theorem \ref{theo:mix_rate} above can be recovered from Ratner's more general result once one remembers the well known fact that the spectrum of the Casimir operator on $L^2(T^1X)$ and the Laplace--Beltrami operator on $L^2(X)$ coincide when restricted to the interval $(-1/4,0)$; beyond this range the spectral data cannot be directly related because of the existence of discrete series.

\section{The Selberg Trace Formula}

\subsection{Motivation.} As we saw in Theorem \ref{theo:mix_rate}, the speed of chaos of the geodesic flow of a closed hyperbolic surface is controlled by the spectral gap of its Laplace--Beltrami operator. But how to compute, or at least estimate, this spectral gap? This question is a rather complicated one, to the point that, even in some of the better studied cases, we do not have a complete answer; see for instance Selberg's famous $1/4$ conjecture \cite{selberg}.

\subsection{The trace formula.} The moral of the story is that accessing the spectral data of the Laplace--Beltrami operator of a closed hyperbolic surface $X$ is quite a formidable task. Perhaps there is hope if one could relate such data to information about the geometry of the surface, which should, at least in spirit, be more accessible to a direct analysis.

Here is where \textit{Selberg's trace formula} enters the picture. Roughly speaking, one could define a natural operator on $L^2(X)$, for instance, using the structure of $X$ as a homogeneous space, and compute its trace in two ways: on one hand, using the geometry of the surface, more precisely, the fact that its closed geodesics are in one-to-one correspondence with conjugacy classes of its fundamental group, and, on the other hand, using the fact that eigenfunctions of the Laplace--Beltrami operator $\Delta$ provide a complete orthonormal system of $L^2(X)$. 
To give a precise statement of this formula we first set up some notation.

Recall that, given a smooth, compactly supported function $H \colon \mathbb{R} \to \mathbb{R}$, its \textit{Fourier transform} is the function $\widehat{H} \colon \mathbb{C} \to \mathbb{C}$ defined as
\[
\widehat{H}(r) := \int_\mathbb{R}H(\ell) \thinspace e^{- i r \ell} \thinspace d \ell.
\]

Given a hyperbolic surface $X$, we will denote by $\mathcal{G}(X)$ the set of \textit{primitive} oriented closed geodesics on $X$; by primitive here we mean closed geodesics that cannot be obtained by concatenating multiple copies of a shorter closed geodesic. Given $\gamma \in \mathcal{G}(X)$, we will denote its hyperbolic length by $\ell(\gamma) > 0$.

We are now ready to state Selberg's trace formula. Although this formula holds in much more generality, we restrict ourselves to the case of closed hyperbolic surfaces.

\begin{theorem}[\cite{sel}]
    \label{theo:selb}
    Let $X$ be a closed hyperbolic surface of genus $g \geq 2$ whose Laplace--Beltrami operator $\Delta$ on $L^2(X)$ has eigenvalues
    \[
    0 = \lambda_0(X) > \lambda_1(X) \geq \lambda_2(X) \geq \cdots \text{ with } \lim_{n \to \infty} \lambda_n(X) = -\infty.
    \]
    For every $j \in \mathbb{N} \cup \{0\}$ denote
    \[
    r_j(X) := \sqrt{-\lambda_j(X) - 1/4} \quad \text{with} \ \  -\pi/2 < \arg(r_j(X)) \leq \pi/2.
    \]
    Then, for every smooth, compactly supported function $H \colon \mathbb{R} \to \mathbb{R}$,
    \[
    \sum_{j = 0}^{+\infty} \widehat{H} (r_j(X)) = (g-1) \int_\mathbb{R} \widehat{H}(r) \tanh(\pi r)  r \thinspace dr + \sum_{\gamma \in \mathcal{G}(X)} \sum_{k=1}^{+\infty} \frac{\ell(\gamma) H(k \ell(\gamma))}{2 \sinh(k \ell(\gamma)/2)}.
    \]
\end{theorem}

For obvious reasons, the left hand side of Selberg's trace formula is usually known as the \textit{spectral side}, while the right hand side is commonly known as the \textit{geometric side}.

\begin{remark}
    The reader interested in an extended discussion on Selberg's trace formula, its derivation, and its consequences, is encouraged to consult \cite{marklof}.
\end{remark}

\subsection{Weyl's law.} Let us take a small detour to illustrate the kind of spectral information one can extract using Selberg's trace formula. Applying a suitable generalization of Theorem \ref{theo:selb} to $\smash{\widehat{H}(r)} = \smash{e^{-\beta r^2}}$ with $\beta \to 0$, one can deduce the following \textit{Weyl law} for the eigenvalues of the Laplace--Beltrami operator of a closed hyperbolic surface:

\begin{corollary}
    Let $X$ be a closed hyperbolic surface of genus $g \geq 2$ whose Laplace--Beltrami operator $\Delta$ on $L^2(X)$ has eigenvalues
    \[
    0 = \lambda_0(X) > \lambda_1(X) \geq \lambda_2(X) \geq \cdots \text{ with } \lim_{n \to \infty} \lambda_n(X) = -\infty.
    \]
    Then,
    \[
    \lim_{\lambda \to +\infty}\frac{\#\{j \in \mathbb{N} \cup \{0\} \thinspace \colon -\lambda_j(X) \leq \lambda\} }{\lambda} = \frac{\mathrm{Area}(X)}{4\pi} = g-1.
    \]
\end{corollary}

\section{Random hyperbolic surfaces}

\subsection{Motivation.} All closed hyperbolic surfaces of a given genus fit together nicely into a \textit{moduli space}. A model of random closed hyperbolic surfaces of a given genus is nothing more than a probability measure in the corresponding moduli space. There are certainly many interesting models one can consider, but for the purposes of this discussion we will restrict ourselves to the \textit{Weil--Petersson} model. Although we will not emphasize it here, this is the natural model to consider from the point of view of hyperbolic geometry.

\begin{remark}
    In parallel to the story we describe in this paper, an equally beautiful study of the spectral properties of random covers of non-compact finite area hyperbolic surfaces was carried out by Hide and Magee \cite{magee}.
\end{remark}

\subsection{The Weil--Petersson model.} To describe the Weil--Petersson model for random closed hyperbolic surfaces of a given genus, we first introduce a local system of coordinates for the corresponding moduli space.

\begin{figure}[ht]
\centering
\includegraphics[scale=.35]{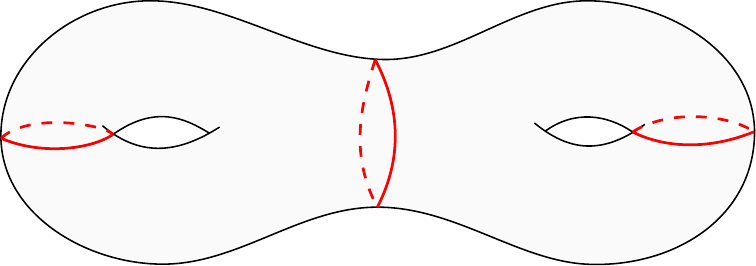}
\caption{A pair of pants decomposition of a genus $2$ surface.}
\label{fig:pdec}
\end{figure}

Notice that, from the point of view of topology, any closed surface of genus $g \geq 2$ can be constructed by gluing $2g-2$ \textit{pairs of pants}, that is, spheres with three boundary components, along their boundaries; see Figure \ref{fig:pdec} for an example of such a decomposition. A similar construction can also be considered for closed hyperbolic surfaces. Indeed, cutting a closed hyperbolic surface of genus $g \geq 2$ along any collection of $3g-3$ disjoint simple closed geodesics, i.e., closed geodesics without self intersections, yields $2g-2$ hyperbolic pairs of pants with geodesic boundary components. This class of pairs of pants is very rigid, in the sense that their geometry is completely determined by the lengths of their boundary components.

\begin{figure}[ht]
		\centering
		\begin{subfigure}{.45\textwidth}
			\centering
			\includegraphics[scale=.35]{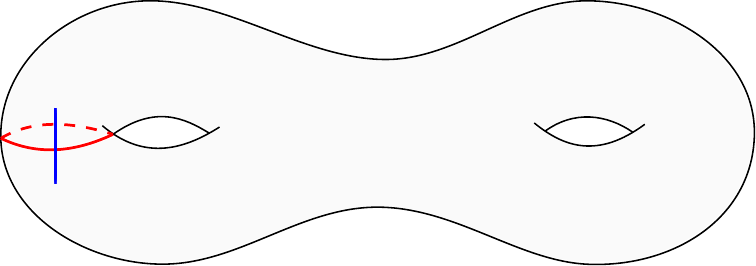}
		\end{subfigure}
		\begin{subfigure}{.45\textwidth}
			\centering
			\includegraphics[scale=.35]{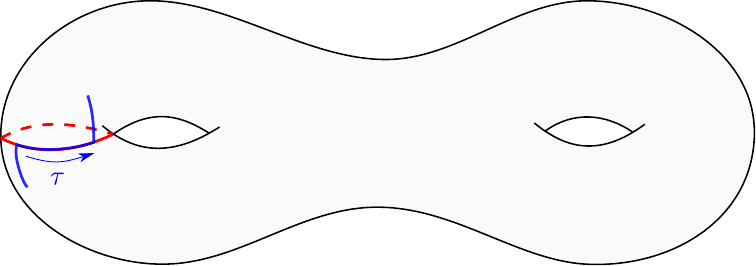}
		\end{subfigure}
        \caption{Twisting a hyperbolic surface along a closed geodesic.}
        \label{fig:fn_tw}
	\end{figure}

Gluing back these pairs of pants following the same original pattern allows us to recover the closed hyperbolic surface we started with. But one needs to be careful at this point, as there are several ways in which one can glue back these pants along their cuffs. Indeed, for each geodesic one cut the original surface along, there is a full circle worth of different twist with which one can glue back the adjacent pairs of pants; see Figure \ref{fig:fn_tw}.

In this way we can provide a local coordinate system near any surface $X$ in the moduli space $\mathcal{M}_g$ of closed hyperbolic surfaces of genus $g \geq 2$. Indeed, given any geodesic pair of pants decomposition of $X$, a local coordinate system is obtained by keeping track of the lengths and the twists, measured in units of hyperbolic length, along the geodesics in the given pair of pants decomposition. These coordinate systems are commonly known as \textit{Fenchel--Nielsen coordinates}. We can now define the desired model of random closed hyperbolic surfaces by appealing to a result of Wolpert:

\begin{theorem}[\cite{wolp}]
    For any $g \geq 2$ there exists a unique finite measure $\mu^{\text{WP}}_g$ on the moduli space $\mathcal{M}_g$, called the \textit{Weil--Petersson measure}, that is locally equal to the Lebesgue measure on any set of Fenchel--Nielsen coordinates.
\end{theorem}

The Weil--Petersson model of random closed hyperbolic surfaces of genus $g \geq 2$ is obtained by normalizing the Weil--Petersson measure $\smash{\mu^{\text{WP}}_g}$ on $\mathcal{M}_g$ to get a probability measure $\smash{\mathbb{P}_g^{\text{WP}}}$. We will denote integrals with respect to this probability measure by $\smash{\mathbb{E}_g^{\text{WP}}}$. It will also be useful to record the total mass of $\smash{\mu^{\text{WP}}_g}$ as the number $V_g > 0$.

\begin{remark}
    The reader interested in the Weil--Petersson geometry of moduli spaces of hyperbolic surfaces and some of its applications is invited to consult \cite{wol_sur}.
\end{remark}

\subsection{Mirzakhani's integration formulas.} Now notice that if we cut a closed hyperbolic surface along a simple closed geodesic we obtain a compact hyperbolic surface with geodesic boundary, potentially disconnected, but certainly of lower complexity. Furthermore, the Weil--Petersson measures on the corresponding moduli spaces behave nicely under this operation because they are explicitely described in terms of length and twist parameters. Using this simple yet powerful inductive idea, Mirzakhani was able to prove remarkable formulas for computing integrals with respect to Weil--Petersson measures.

Before stating a particular instance of these formulas, let us first introduce some notation. Given a closed hyperbolic surface $X$, denote by $\mathcal{G}_s(X) \subseteq \mathcal{G}(X)$ the set of all oriented simple closed geodesics on $X$. Given a measurable function $F \colon \mathbb{R} \to \mathbb{R}$ with compact support, consider the sum
\[
F_s(X) := \sum_{\gamma \in \mathcal{G}_s(X)} F(\ell(\gamma)).
\]

The following is just one of the many variants of Mirzakhani's integration formulas:

\begin{theorem}[\cite{mir_ct}]
    \label{theo:mir_ct}
    For every $g \geq 2$ there exists an explicit polynomial function $V_g^s \colon \mathbb{R} \to \mathbb{R}$ of degree $6g-7$, called a \textit{volume polynomial}, with the following property. Let $F \colon \mathbb{R} \to \mathbb{R}$ be a continuous function with compact support. Then,
    \[
    \mathbb{E}_g^{\text{WP}}\left(F_s(X)\right) = \frac{1}{V_g} \int_0^{+\infty}F(\ell) V_g^s(\ell) \thinspace d\ell.
    \]
\end{theorem}

\begin{remark}
    Although we will not elaborate further here, it is a crucial part of Mirzakhani's work that the functions $V_g^s(\ell)$ can be computed recursively by induction on the integer $g$.
\end{remark}

Notice how fascinating Theorem \ref{theo:mir_ct} is: We have a very complicated moduli space of closed hyperbolic surfaces, yet, somehow, we are able to compute integrals on this space via integration with respect to an explicit polynomial density on $\mathbb{R}$!

To illustrate how explicit the description of the volume polynomials $V_g^s(\ell)$ in Theorem \ref{theo:mir_ct} can be, we present the following estimate of Mirzakhani and Petri:

\begin{theorem} [\cite{mirpet}]
    \label{theo:mir_pet}
    There exists $c > 0$ such that uniformly for all $g \geq 2$ and all $\ell > 0$,
    \[
    \frac{V_g^s(\ell)}{V_g} = \frac{4}{\ell} \sinh^2\left( \frac{\ell}{2}\right) + O\left( \frac{(1+\ell)^c e^\ell}{g}\right).
    \]
\end{theorem}

\subsection{Mirzkhani's counting formulas.} Although we will not explore these ideas any further, we would like to at least mention in passing what is perhaps one of the most famous applications of Mirzakhani's integration formulas in Theorem \ref{theo:mir_ct}, that is, her counting formulas for simple closed geodesics on closed hyperbolic surfaces:

\begin{theorem} [\cite{mir_ct}]
    \label{theo:ct}
    For every closed hyperbolic surface $X$ of genus $g \geq 2$ there exists a positive constant $B(X) > 0$ such that
    \[
    \lim_{L \to \infty} \frac{\#\{\gamma \in \mathcal{G}_s(X) \thinspace \colon \thinspace \ell(\gamma) \leq L\}}{L^{6g-6}} = B(X).
    \]
\end{theorem}

\begin{remark}
    The reader interested in exploring more about the history of counting results in the direction of Theorem \ref{theo:ct} is invited to consult \cite{surv}.
\end{remark}

\section{Optimal spectral gap}

\subsection{Main result.} We are now ready to state a precise version of the following result highlighted in \S \ref{sec:intro}: Random large genus closed hyperbolic surfaces have optimal spectral gap. Notice that, by Theorem \ref{theo:opt}, one cannot hope for spectral gaps of large genus closed hyperbolic surfaces to be greater than $1/4$. In this context, the main result of Anantharaman and Monk we wish to discuss is the following:

\begin{theorem} [\cite{gap1,monksur,gap2}]
    \label{theo:main}
    For every $\epsilon > 0$,
    \[
    \lim_{g \to \infty} \mathbb{P}_g^{\text{WP}}\left( -\lambda_1(X) \leq \frac{1}{4} - \epsilon\right) = 0.
    \]
\end{theorem}

\begin{remark}
    An analogous result for random regular graphs of a fixed valency with many vertices was conjectured first by Alon \cite{alon} and later proved by Friedman \cite{friedman}. Theorem \ref{theo:main} and parts of its proof are directly inspired by this result.
\end{remark}

\subsection{Towards optimality.} Let us describe a basic strategy one can try to follow to prove Theorem \ref{theo:main} using the tools introduced in previous sections.

\begin{figure}[ht]
\centering
\begin{subfigure}{.45\textwidth}
\centering
\includegraphics[scale=.18]{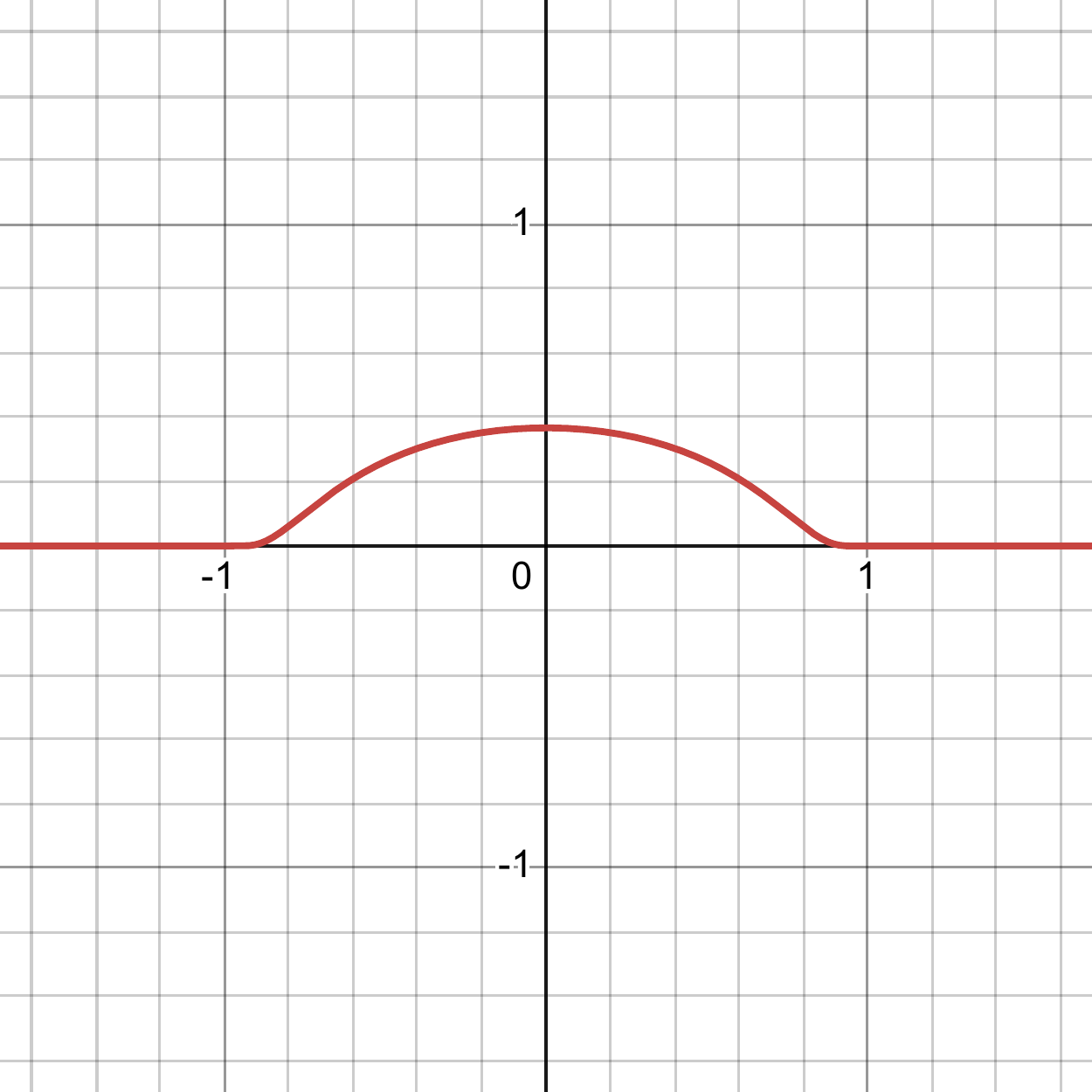}
\end{subfigure}
\begin{subfigure}{.45\textwidth}
\centering
\includegraphics[scale=.133]{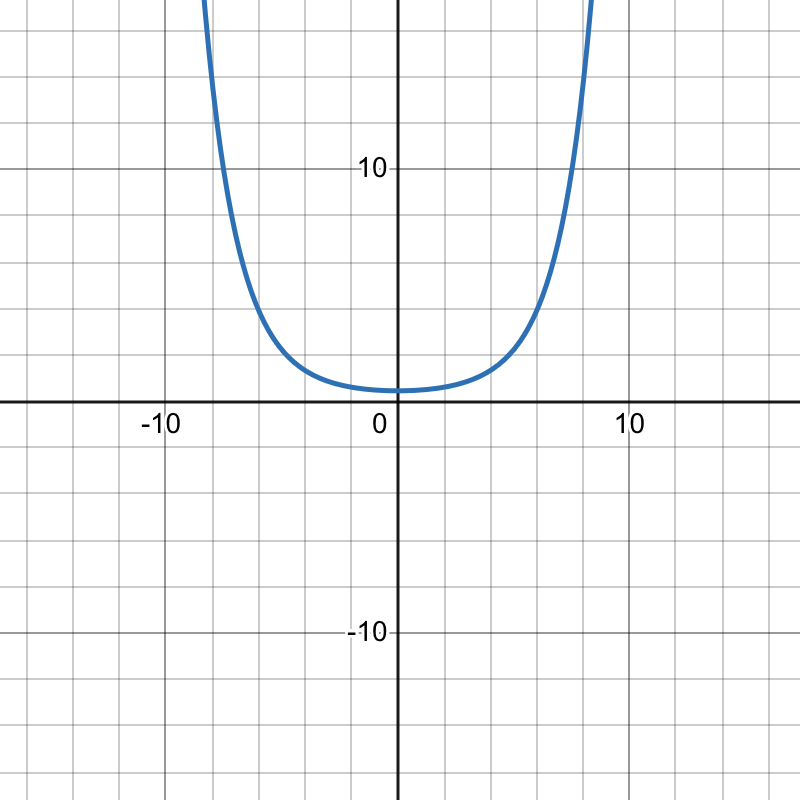}
\end{subfigure}
\caption{The test function $H$ (left) and the restriction of its Fourier transform $\widehat{H}$ to the imaginary axis (right).}
\label{fig:test}
\end{figure}

Our starting point will be Selberg's trace formula; see Theorem \ref{theo:selb} for a precise statement. For the moment, fix an arbitrary closed hyperbolic surface $X$ of genus $g$; later this surface will be chosen at random. Consider a test function $H \colon \mathbb{R} \to \mathbb{R}$ that is smooth, non-negative, even, and with support equal to $[-1,1]$; see Figure \ref{fig:test}. Furthermore, assume its Fourier transform $\smash{\widehat{H}} \colon \mathbb{C} \to \mathbb{C}$ is real valued and non-negative on $\mathbb{R} \cup i [-1/2,1/2]$. For large values of $L \geq 1$ we will apply the trace formula using the test functions
\[
H_L(\ell) := H(\ell/L) \quad \text{and} \quad \widehat{H_L}(r) = L\widehat{H}(rL).
\]
Notice that $H_L$ is supported on $[-L,L]$, so the geometric side of Selberg's trace formula will only involved closed geodesics $\gamma \in \mathcal{G}(X)$ with $\ell(\gamma) \leq L.$

The first key observation is that eigenvalues $\lambda_j > -1/4$ of the Laplace--Beltrami operator $\Delta$ on $L^2(X)$ correspond to purely imaginary values $r_j \in i[-1/2,1/2] \setminus \{0\}$ in Theorem \ref{theo:selb}. In particular, because $H$ is non-negative, the corresponding terms on the spectral side of Selberg's trace formula can be bounded roughly as
\[
\widehat{H_L}(r_j) = 2L \int_{0}^{1} H\left(\ell \right) \thinspace e^{|r_j| \ell} \thinspace d\ell\geq \left( \int_{1/2}^{1} H(\ell) \thinspace d\ell  \right) \cdot e^{|r_j|L/2}.
\]
In other words, because $H$ has support equal to $[-1,1]$, small eigenvalues of the Laplace--Beltrami operator $\Delta$ on $L^2(X)$, whose existence we are trying to rule out, correspond to exponential growth on the spectral side of Selberg's trace formula. Thus, at least intuitively, we should aim to show that the geometric side of the trace formula grows subexponentially.

Denote $r_1$ by $r_1(X)$, making the dependence on $X$ explicit. Being a bit more precise with the bounds above, one can show that, for every $\alpha > 0$ and every $\epsilon > 0$, there exists a constant $C_{\alpha,\epsilon} > 0$ such that for every closed hyperbolic surface $X$,
\[
-\lambda_1(X) \leq  \frac{1}{4} - \alpha^2 - \epsilon \quad \Rightarrow \quad \widehat{H_L}(r_1(X)) \geq C_{\alpha,\epsilon} \thinspace e^{(\alpha+\epsilon) L}.
\]
In particular, for every $g \geq 2$,
\[
\mathbb{P}_{g}^\mathrm{WP}\left( -\lambda_1(X) \leq  \frac{1}{4} - \alpha^2 - \epsilon\right) \leq  \mathbb{P}_{g}^\mathrm{WP} \left( \widehat{H_L}(r_1(X)) \geq C_{\alpha,\epsilon} \thinspace e^{(\alpha+\epsilon) L}\right).
\]
The role of $\alpha$ here is to represent different degrees of precision one can aim for when proving Theorem \ref{theo:main}; the goal is to prove the desired statement for arbitrarily small $\alpha$.

As is usual in probability theory, and recalling the non-negativity assumption on $\smash{\widehat{H}}$, one can use Markov's inequality to bound
\[
\mathbb{P}_{g}^\mathrm{WP} \left( \widehat{H_L}(r_1(X)) \geq C_{\alpha,\epsilon} \thinspace e^{(\alpha+\epsilon) L}\right) \leq \frac{\mathbb{E}_{g}^\mathrm{WP}\left( \widehat{H_L}(r_1(X)) \right)}{C_{\alpha,\epsilon} \thinspace e^{(\alpha+\epsilon)L}}.
\]
This reduces the proof of Theorem \ref{theo:main} to showing that, for some appropriate choice of $L = L(g)$, later to be chosen so that $L(g) \to +\infty$ as $g \to \infty$,
\[
\lim_{g \to \infty} \frac{\mathbb{E}_{g}^\mathrm{WP}\left( \widehat{H_{L(g)}}(r_1(X)) \right)}{C_{\alpha,\epsilon} \thinspace e^{(\alpha+\epsilon)L(g)}} = 0.
\]

To derive these asymptotics on the expected value we integrate the trace formula. But notice first that on both sides of this formula there are terms that do not depend on the underlying closed hyperbolic surface. On the spectral side, this corresponds to the eigenvalue $\lambda_0 = 0$, or, equivalently, to $r_0 = i/2$. Indeed, we have
\begin{equation}
\label{eq:A}
\widehat{H_L}(r_0) =  \widehat{H_L}(i/2) = 2 \int_0^{+\infty} H_L(\ell) \cosh(\ell/2) \thinspace d \ell.
\end{equation}
On the geometric side, the constant term corresponds to the one outside the sum over all closed geodesics. Recall that $\tanh(x) \leq 1$ for all $x \in \mathbb{R}$, and that, because $H$ has compact support, $\smash{\widehat{H}}$ decays faster than any polynomial. Thus, we can bound
\[
\bigg\vert(g-1) \int_\mathbb{R} \widehat{H}_L(r) \tanh(\pi r)  r \thinspace dr \bigg\vert \leq g \|r \widehat{H}_L(r)\|_{L^1(\mathbb{R})} = \frac{g}{L}\|r \widehat{H}(r)\|_{L^1(\mathbb{R})}.
\]

Another important preliminary observation to make is that the terms coming from non-primitive closed geodesics, i.e., from $k \geq 2$ in Theorem \ref{theo:selb}, are of lower order. Although we will not go into details here, an appropriate estimate can be obtained by using the fact that the $\sinh(k \ell(\gamma)/2)$ factor in the denominator of the sum of the geometric side of the trace formula grows exponentially with $k$, and so the corresponding contributions will be much smaller for $k \geq 2$ than for $k =1$. More precisely, one can show there exists a constant $C' > 0$ such that for any $g \geq 2$,
\[
\mathbb{E}_{g}^\mathrm{WP}\left( \sum_{\gamma \in \mathcal{G}(X)} \sum_{k=2}^{+\infty} \frac{\ell(\gamma) H_L(k \ell(\gamma))}{2 \sinh(k \ell(\gamma)/2)} \right) \leq C' \|H\|_\infty L^2 g.
\]

The contribution of simple closed geodesics will also play an important role in our discussion. Using Mirzakhani's integration formulas, i.e., Theorem \ref{theo:mir_ct}, and the estimates on volume polynomials in Theorem \ref{theo:mir_pet}, one can show that, uniformly for all $g \geq 2$,
\[
\mathbb{E}_{g}^\mathrm{WP}\left( \sum_{\gamma \in \mathcal{G}_s(X)} \frac{\ell(\gamma) H_L( \ell(\gamma))}{2 \sinh( \ell(\gamma)/2)} \right) = 2 \int_0^{+\infty} H_L(\ell) \cosh(\ell/2) \thinspace d \ell + O\left(1 + \frac{L^c e^{L/2}}{g} \right).
\]
Notice how the leading term in this estimate, which contributes to the geometric side of the trace formula, exactly matches the spectral side computation in \eqref{eq:A}.

From this discussion it follows that, after integrating Selberg's trace formula  over moduli space, one can deduce that, uniformly on $g \geq 2$,
\[
\mathbb{E}_{g}^\mathrm{WP}\left( \widehat{H_{L}}(r_1(X)) \right) \leq O\left(L^2 g + \frac{L^c e^{L/2}}{g} \right) + \mathbb{E}_{g}^\mathrm{WP}\left( \sum_{\gamma \in \mathcal{G}(X)\ \setminus \mathcal{G}_s(X)} \frac{\ell(\gamma) H_L( \ell(\gamma))}{2 \sinh( \ell(\gamma)/2)} \right).
\]
It turns out that, through very sophisticated geometric arguments, one can show that, for every $\epsilon > 0$, there exists $C_\epsilon'' > 0$ such that, uniformly on $g \geq 2$,
\begin{equation}
\label{eq:B}
\mathbb{E}_{g}^\mathrm{WP}\left( \sum_{\gamma \in \mathcal{G}(X)\ \setminus \mathcal{G}_s(X)} \frac{\ell(\gamma) H_L( \ell(\gamma))}{2 \sinh( \ell(\gamma)/2)} \right) \leq C_\epsilon'' \frac{e^{(1+\epsilon)L/2}}{g}.
\end{equation}
In particular, if we take $\alpha = 1/4$ and  $L(g) = 4 \log g$, we conclude
\[
\lim_{g \to \infty} \frac{\mathbb{E}_{g}^\mathrm{WP}\left( \widehat{H_{L(g)}}(r_1(X)) \right)}{C_{\alpha,\epsilon} \thinspace e^{(\alpha+\epsilon)}} = 0.
\]

This is precisely the approach considered independently by Wu and Xue \cite{gp1} and Lipnowski and Wright \cite{gp2} in their proofs of the weaker version of Theorem \ref{theo:main} with a spectral gap of $3/16$. Let us highlight that their arguments to prove \eqref{eq:B} are quite different. 

\subsection{Optimality.} All that work and we are still far off from the desired spectral gap of $1/4$. What needs to be improved to reach that threshold? On one hand, we cannot expect a bound as in \eqref{eq:B} to be sufficient: we should aim for better estimates on the corresponding sums over non-simple closed geodesics. Furthermore, even the estimates in Theorem \ref{theo:mir_pet} for simple closed geodesics will not be enough for our purposes: we will need more precise expansions in inverse powers of $g$.

\subsubsection{Alternative trace scheme.} 
\begin{figure}[ht]
\centering
\begin{subfigure}{.45\textwidth}
\centering
\includegraphics[scale=.18]{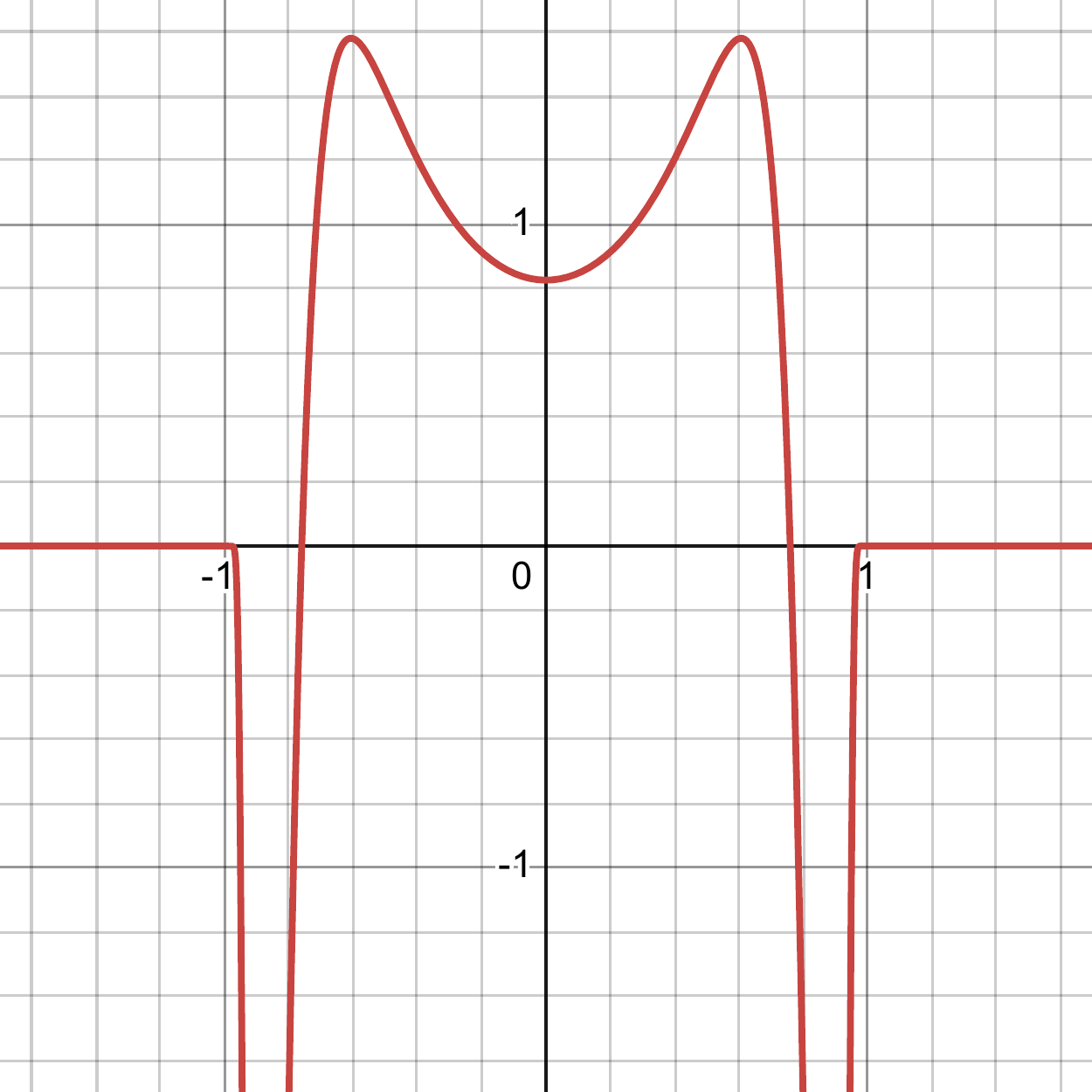}
\end{subfigure}
\begin{subfigure}{.45\textwidth}
\centering
\includegraphics[scale=.133]{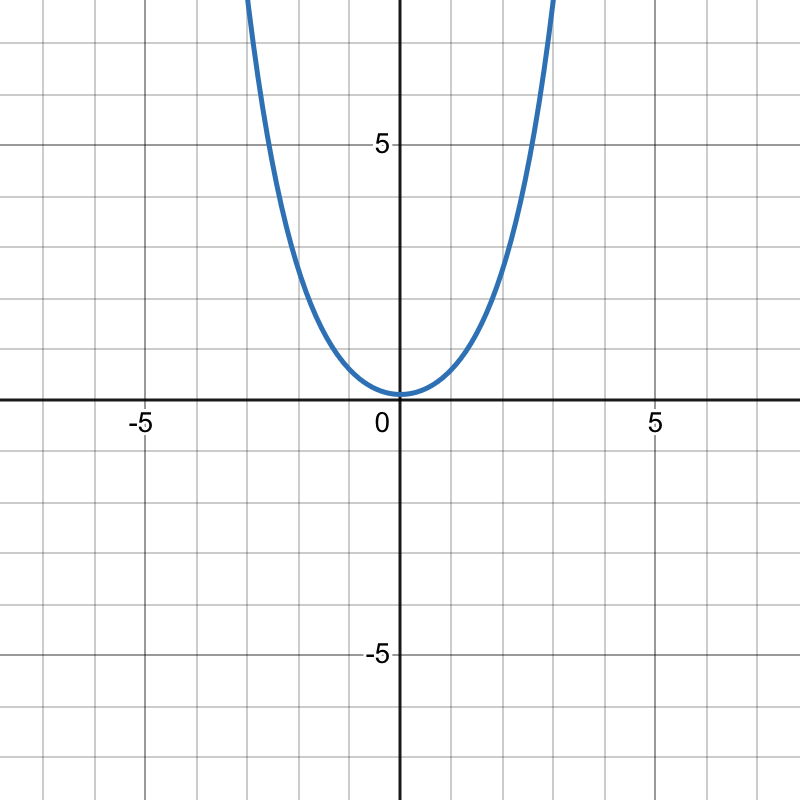}
\end{subfigure}
\caption{The test function $\mathcal{D}H$ (left) and the restriction of its Fourier transform $\widehat{\mathcal{D}H}$ to the imaginary axis (right).}
\label{fig:bump2}
\end{figure}

Before we proceed any further, let us introduce an alternative application of Selberg's trace formula that will be more suitable for our purposes. The main aim is to kill the contributions coming from the trivial eigenvalue $\lambda_0 = 0$, i.e., $r_0 = i/2$.  To do so we consider the differential operator 
\[
\mathcal{D} := \frac{1}{4} - \frac{d^2}{d\ell^2}
\]
and for some arbitrary $m \geq 1$ we replace the test functions above by
\[
\mathcal{D}^m H_L(\ell) \quad \text{and} \quad \widehat{\mathcal{D}^m H_L}(r) = \left( \frac{1}{4} + r^2\right)^m  \widehat{H_L}(r).
\]
This operation does not change the support of the test function under consideration. Nevertheless, although the Fourier transform of the test function still retains the desired non-negativity properties, the test function itself now changes signs; see Figure \ref{fig:bump2}. This will introduce significant difficulties that will need to be dealt with later.

Another important point to make is that, in killing the contributions of the trivial eigenvalue, we produced test functions for which eigenvalues near zero do not longer give rise to uniform exponential growth. Nevertheless, as we already established a uniform spectral gap for random large genus hyperbolic surfaces, this issue can be easily bypassed; in fact, it is enough to quote Mirzakhani's softer original proof of uniform spectral gap \cite{tan1}.

\subsubsection{General integration formulas.} 
\begin{figure}[ht]
\centering
\includegraphics[scale=.35]{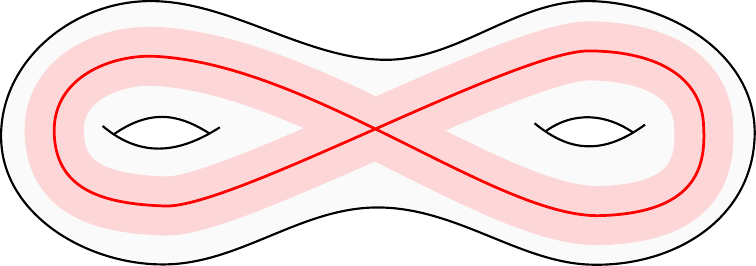}
\caption{The local topological type of a closed geodesic.}
\label{fig:tt}
\end{figure}

Back to the first of our main concerns: How do we even integrate functions of non-simple closed geodesics if Mirzakhani's inductive approach only seems to work for simple ones? Let us introduce the first main insight of Anantharaman and Monk, \textit{local topological types}. Roughly speaking, the local topological type of a closed geodesic on a closed hyperbolic surface is the isotopy class of the subsurface obtained by thickening the geodesic to a tubular neighborhood together with the information of how the geodesic fills such subsurface; see Figure \ref{fig:tt}. The most basic local topological type are simple closed geodesics: they all thicken to an annulus. Notice that we can make sense of a given local topological type on any surface of sufficiently large genus.

Given a closed hyperbolic surface $X$ and a local topological type of closed geodesics $T$, denote by $\mathcal{G}_T(X) \subseteq \mathcal{G}(X)$ the set of all oriented closed geodesics on $X$ of local topological type $T$. Just as for simple closed geodesics, given a continuous function $F \colon \mathbb{R} \to \mathbb{R}$ with compact support, consider the sum
\[
F_T(X) := \sum_{\gamma \in \mathcal{G}_T(X)} F(\ell(\gamma)).
\]

Although general closed geodesics might not be simple, and thus, a priori, not amenable with inductive procedures, the geodesic representatives of the boundary components of the surfaces obtained via thickening are always simple and thus allow one to prove the following result in the spirit of Mirzakhani's integration formulas:

\begin{theorem}[\cite{gap1}]
    \label{theo:am_ct}
    Let $T$ be a local topological type of closed geodesics. Then, for every $g \geq 2$ there exists a function $V_g^T \colon \mathbb{R} \to \mathbb{R}$,  called a \textit{volume function}, with the following property. Let $F \colon \mathbb{R} \to \mathbb{R}$ be a measurable function with compact support. Then,
    \[  \mathbb{E}_g^{\text{WP}}\left(F_T(X)\right) = \frac{1}{V_g} \int_0^{+\infty}F(\ell) V_g^T(\ell) \thinspace d\ell.
    \]
\end{theorem}

Notice that, in contrast with Theorem \ref{theo:mir_ct}, the volume function $V_g^T(\ell)$ in Theorem \ref{theo:am_ct} need not be a polynomial. Nevertheless, through sophisticated geometric methods, one can still control its behavior in a very precise way, as we now discuss.

\subsubsection{Friedman-Ramanujan functions.} We would like to have an estimate as in Theorem \ref{theo:mir_pet} for volume functions $V_g^T(\ell)$ of general local topological types $T$. Furthermore, we would like an even more precise expansion in inverse powers of $g$. Relying on previous results of Mirzakhani and Zograf \cite{ZoMir}, such an expansion was proven to exist in the work of Anantharaman and Monk:

\begin{theorem}[\cite{gap1}]
    \label{theo:expansion}
     Let $T$ be a local topological type of closed geodesics corresponding to a filled surface $S$ of Euler characteristic $\chi(S) \leq 0$. Then, the volume function $V_g^T(\ell)$ admits an asymptotic expansion in inverse powers of $g$ in the following sense: There exists a unique family of continuous functions $(f^T_k \colon \mathbb{R} \to \mathbb{R})_{k \geq -\chi(S)}$ such that for any integer $K \geq 0$, any large enough $g \geq 2$, and any $\eta > 0$,
     \[
     \frac{V_g^T(\ell)}{V_g} = \sum_{k = -\chi(S)}^K \frac{f_k^T(\ell)}{g^k} + O_{K,S,\eta} \left( \frac{e^{(1+\eta)\ell}}{g^{K+1}}\right).
     \]
\end{theorem}

This expansion would not be of much use without a good understanding of the structure of the functions $f^T_k(\ell)$. Inspired by Friedman's proof of Alon's conjecture \cite{friedman}, Anantharaman and Monk introduced a class of so-called Friedman--Ramanujan functions:

\begin{definition}
\label{def:class}
A continuous function $r \colon \mathbb{R}_{>0}\to \mathbb{C}$ is said to be of class $\mathcal{R}_w$ if there exist constants $c_1,c_2 > 0$  such that for all $L \geq 1$,
\[
\int_0^L |r(s)| \thinspace ds \leq c_1 (n+1)^{c_2} e^{L/2}.
\]
A continuous function $f \colon \mathbb{R}_{>0}\to \mathbb{C}$ is said to be a \textit{Friedman--Ramanujan function} if there exists a polynomial function $p \colon \mathbb{R}\to \mathbb{C}$ and a function $r \colon \mathbb{R}_{>0} \to \mathbb{C}$ of class $\mathcal{R}_w$ such that
\[
f(\ell) = p(\ell) e^\ell + r(\ell), \quad \text{for all } \ell > 0.
\]
\end{definition}

Roughly speaking, Friedman--Ramanujan functions produce extra cancellations when integrated against the $\sinh(\ell/2)$ denominator in Selberg's trace formula when using the alternative scheme proposed above. Such cancellations arise via integration by parts once one notices that, for any such function, if the leading term corresponds to a polynomial of degree strictly less than $m$, then the 
product of such leading term with $\smash{e^{-\ell/2}}$ belongs to the kernel of the differential operator $\smash{\mathcal{D}^m}$.

Perhaps the most important and technical step in the proof of Anantharaman and Monk of Theorem \ref{theo:main} is verifying that the functions $f_k^T(\ell)$ in the expansions provided by Theorem \ref{theo:expansion} are Friedman--Ramanujan functions:

\begin{theorem}[\cite{gap2}]
    \label{theo:FR}
    Let $T$ be a local topological type of closed geodesics corresponding to a filled surface $S$ of Euler characteristic $\chi(S) \leq 0$. Define $s := 1$ if $T$ corresponds to simple closed geodesics and $s := 0$ otherwise. Then, for any $k \geq -\chi(S)$, the function 
    \[
    \ell \in \mathbb{R}_{>0} \mapsto\ell^s f_k^T(\ell) \in \mathbb{R},
    \]
    with $f_k^T(\ell)$ as in Theorem \ref{theo:expansion}, is a Friedman--Ramanujan function.
\end{theorem}

\begin{remark}
    A more precise version of Theorem \ref{theo:FR} controlling the degree of the corresponding polynomial and the constants $c_1,c_2 > 0$ in Definition \ref{def:class} in terms of $k$ and $T$ is needed and proved in the work of Anantharaman and Monk. With the aim of simplifying the exposition here, we have opted to restrict our attention to the weaker version above.
\end{remark}

But how to show that a function is Friedman--Ramanujan? It turns out that this class of functions can be characterized as the solution set of a simple integro-differential equation. Using this characterization, one can show stability of this class under certain convolution operations with respect to a class of sufficiently regular kernels. This class of kernels contains variants of hyperbolic trigonometric functions that naturally show up when computing lengths of closed geodesics on compact hyperbolic surfaces. The proof of Anantharaman and Monk involves constructing variants of Fenchel--Nielsen coordinates of compact hyperbolic surfaces for which the lengths of filling closed geodesics can be represented as convolutions with respect to the right class of kernels.

\subsubsection{Tangles.} We seem to be on the right track, but there is a glaring issue we have not yet addressed. Suppose that we proceed with our alternative trace scheme, and, inspired by the suboptimal gap bound above, we consider $L(g) := A \log g$ for some large but fixed constant $A > 0$. One can show that, at this scale, when analizing the geometric side of Selberg's trace formula, it is enough to consider closed geodesics filling subsurfaces of Euler characteristic at least $-(2A + 1)$. Nevertheless, even with this bound on the Euler characteristic of the filled subsurface, there could still be infinitely many local topological types to consider, all contributing error terms that will not be summable in the end.

\begin{figure}[ht]
\centering
\includegraphics[scale=.5]{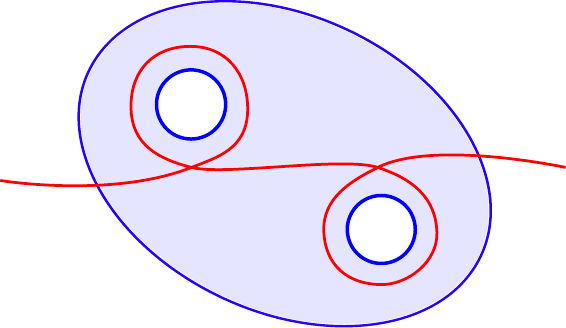}
\caption{A tightly looping geodesic (red) and its tangle (blue).}
\label{fig:tangle}
\end{figure}

What leads to the proliferation of too many local topological types? This phenomena is explained by Anantharaman and Monk via the notion of \textit{tangles}; see Figure \ref{fig:tangle}. According to this picture, a closed geodesic can wind around a surface in tight loops in such a way that it gives rise to a geodesic pair of pants with short boundary components. More generally, a tangle is an extremely short closed geodesic, a geodesic pair of pants with short boundary components, or a geodesic once-holed tori with short boundary.

Denote by $\mathrm{TF}_g \subseteq \mathcal{M}_g$ the set of tangle-free closed hyperbolic surfaces of genus $g \geq 2$. Technically speaking, a proper definition of tangle-free surfaces should involve parameters that control the shortness of the geodesics in the bad configurations specified above; we will not aim for such level of precision in our discussion. From work of Mirzakhani, Monk, and Thomas, one can deduce that random large genus closed hyperbolic surfaces are tangle-free in the following sense:

\begin{theorem} \cite{tan1,tan2}
    For an appropriate choice of parameters in the definition of $TF_g$, the following limit can be made arbitrarily small with an explicit convergence rate, 
    \[
    \lim_{g \to \infty} \mathbb{P}_g^{\mathrm{WP}}\left( \mathrm{TF}_g \right).
    \]
\end{theorem}

Furthermore, as one would hope from our previous discussion, the proliferation of local topological types in tangle-free surfaces can be controlled:

\begin{theorem} \cite{moeb}
    Fix a large constant $A > 0$ and $g \geq 2$. Then, the number of local topological types that can be represented as a filling closed geodesic $\gamma$ of length $\ell(\gamma) \leq A \log g$ on a tangle-free compact hyperbolic surface of Euler characteristic at least $-(2A+1)$ is bounded uniformly by a polynomial in $\log g$.
\end{theorem}

\subsubsection{Moebius inversion.} Everything seems to be in place for the proof of Theorem \ref{theo:main} to go through. But in handling the proliferation of too many local topological types we have introduced a new issue: the indicator function $\mathbbm{1}_{\mathrm{TF}_g}$ cannot be integrated a la Mirzakhani. Notice also that we cannot just use the crude bound $\mathbbm{1}_{\mathrm{TF}_g} \leq 1$ because the test function considered in the alternative trace scheme changes signs.

Do not despair though, as a well known trick also featured in the work of Friedman on random graphs can be used to address this issue. The idea is to perform inclusion-exclusion to rewrite the indicator function $\mathbbm{1}_{\mathrm{TF}_g}$ as an alternating sum of geometric counting functions. In general, if $N$ is a random variable counting the number of a appearances of a geometric pattern, in our case $N$ will count tangles, then one can rewrite
\begin{equation}
\label{eq:moeb}
\mathbbm{1}_{\{N = 0\}} = 1 - \sum_{j=1}^{+\infty} \frac{(-1)^{j+1}}{j!} N_j,
\end{equation}
where, for every $j \in \mathbb{N}$, the random variable $N_j$ corresponds to the number of ordered families of $j$ geometric patterns that are being counted by $N$. Notice these random variables are counting functions, and so one hopes to be able to integrate them as above.

Nevertheless, tangles are relatively complicated geometric patterns, in the sense that the variables $N_j$ do not have a clear geometric interpretation. Hence, one cannot expect to apply a formula as clean as \eqref{eq:moeb} directly in our setting. Nevertheless, Anantharaman and Monk are able to prove a similar, albeit less explicit, \textit{Moebius inversion formula} that can be integrated in the sense of Mirzakhani. This formula is quite involved and describing it explicitly would escape the scope of our discussion; see \cite{moeb} for a precise statement.

\subsubsection{Finishing the proof.} All the ingredients are now in place to finish the proof of Theorem \ref{theo:main}. After discarding surfaces with tangles using the Moebius inversion formula, the polynomially many relevant terms on the geometric side of Selberg's trace formula can be integrated to obtain cancelations with summable errors. More precisely, given an arbitrary integer $K \geq 0$, one considers the scale $L(g) := 2(K+2) \log g$ and fixes an integer $m$ depending on the degrees of the polynomials in the leading terms of the Friedman--Ramanujan functions that show up in the expansion of order $K$ provided by Theorem \ref{theo:expansion}. This way one concludes that for $\alpha := 1/2(K+2)$ and every $\epsilon > 0$,
    \[
    \lim_{g \to \infty} \mathbb{P}_g^{\text{WP}}\left( -\lambda_1(X) \leq \frac{1}{4} - \alpha^2 -\epsilon\right) = 0.
    \]

\subsection{Conclusion.} In closing, we would like to invite the reader to look back at the main problem at hand, its history, and its implications. We are interested in estimating the spectral gap of the Laplace--Beltrami operator of a typical closed hyperbolic surface of large genus, with the hope of showing that it is arbitrarily close to the optimal bound of $1/4$. This is relevant because, among many other things, such spectral gap controls the speed of chaos of the corresponding geodesic flow.

But experience shows that spectral data is very hard to access, even in explicit examples. Regardless, we can keep our hopes high and try to inspire ourselves on Selberg's success story. This leads us to use his trace formula to relate spectral data to geometric data that should, at least a priori, be more easily accessible. As we only care about typical  closed hyperbolic surfaces, using any type of Markov inequality we can reduce our problem to estimating the integral of the trace formula over the corresponding moduli space.

Now, although originally introduced to address a completely different set of problems, we can use Mirzakhani's integration formulas to compute the integral of the geometric side of the trace formula. But only when the closed geodesics involved are simple, or so one would expect at first glance. Not too bad, as this lead to a spectral gap of $3/16$.

Enter now the revolution of Anantharaman and Monk, who teach us that, although explicit formulas à la Mirzakhani a priori do not seem attainable for moduli space integrals involving non-simple closed geodesics, it is actually possible to compute precise estimates by appealing to the notion of Friedman--Ramanujan functions. Furthermore, the cancellations obtained via this approach lead to summable error terms as long as one discards closed hyperbolic surfaces with tangles. Although not possible in the most naive way, such surfaces can be discarded via a sophisticated Moebius inversion formula. All of this can be carried out with such precision as to obtain the desired optimal spectral gap of $1/4$.

Before going to bed tonight, just contemplate for a moment the beauty of all of these ideas from seemingly unrelated sources, geometry, analysis, and probability, coming together to provide a complete and clean answer to such an ellusive and difficult problem. And do not forget to reflect on the creativity and technical prowess that was needed by Anantharaman and Monk to see this \textit{tour de force} proof to the end!

\bibliographystyle{amsalpha}


\bibliography{bibliography}

$ $

\end{document}